\documentclass[12pt, letterpaper]{article}
\usepackage{amssymb}
\usepackage{amsfonts}
\usepackage{amsmath,amsthm}
\usepackage[unicode]{hyperref}
\usepackage{tikz}
\usepackage{graphicx}
\usepackage{xcolor}

\newtheorem{thm}{Theorem}[section]
\newtheorem{prop}[thm]{Proposition}
\newtheorem{lem}[thm]{Lemma}
\newtheorem{rem}[thm]{Remark}
\newtheorem{Cor}[thm]{Corollary}

\makeindex

\title{Normal subgroups in the  group of column-finite infinite matrices}

\author{Waldemar Ho{\l}ubowski\footnote{Corresponding author, Email address: w.holubowski@polsl.pl}, Martyna Maciaszczyk\\ and Sebastian \.Zurek}

\date{}

\begin{document}

\maketitle

Silesian University of Technology, Faculty of Applied Mathematics\\
Kaszubska 23, 44-101 Gliwice, Poland

\begin{abstract}
The classical result, due to Jordan, Burnside, Dickson, says that every normal subgroup of $GL(n, K)$ ($K$ - a field, $n \geq 3$) which is not contained in the center, contains $SL(n, K)$. A. Rosenberg gave description of normal subgroups of $GL(V)$, where $V$ is a vector space of any infinite cardinality dimension over a division ring.  However, when he considers subgroups of the direct product of the center and the group of linear transformations $g$ such that $g-id_V$ has  finite dimensional range the proof is not complete. We fill this gap  for countably dimensional $V$ giving description of the lattice of normal subgroups in the group of infinite column-finite matrices indexed by positive integers over any field.
\end{abstract}



\section{Introduction}
Description of normal subgroups is a fundamental problem in group theory. The classical result due to Jordan, Burnside, Dickson, says that every normal subgroup of $GL(n, K)$ ($K$ - a field, $n \geq 3$) which is not contained in the center (the group of scalar matrices $D_{sc}(n, K)$), contains $SL(n, K)$. Moreover, the factorgroup $SL(n, K)/ (SL(n, K) \cap D_{sc}(n, K))$ is simple. The last groups, in case of finite field $K= \mathbb{F}_q$, define an infinite family of finite simple groups 
$PSL(n, \mathbb{F}_q)$. 
 We note that  $SL(n,K) \cap D_{sc}(n, K)$ can be nontrivial and depends on $n$ and characteristic of $K$.
For example, if $K=\mathbb{F}_p$ ($p-$prime) then the scalar matrix ${\rm diag} (x, x, \ldots , x)$ for nonzero $x$ belongs to $SL(p-1, \mathbb{F}_p)$ (we have $x^{p-1}=1$).

In \cite{Dieu} J. Dieudonn\'e extended above result to the case $GL(n, D)$, where $n \geq 3$ and $D$ is any division ring. In fact, this paper  defines noncommutative determinant for $GL(n, D)$ and study the properties of its commutator subgroup. The last paragraph of this paper (number 14 on page 45)  announce similar results for infinite dimensional case, contains only three sentences and no proofs. 

Let $V$ be a linear space of infinite dimension $\aleph_{\delta}$ over a division ring $D$. By $End(V)$ we denote the ring of endomorphisms of $V$ and by $GL(V)$ its group of units, i.e. the group of all invertible linear transformations on $V$.

Dieudonn\'e consider the subgroup $C(V)$ of $GL(V)$ which is generated by transvections, i.e. elementary transvections and their conjugates,  and says that $C(V)$ coincides with commutator subgroup of $F(V)$ -- the group of linear transformations $g$ such that $g-id_V$ has  finite dimensional range.  Normal subgroups of $F(V)$ which are not contained in  the center $Z(V)$ contain $C(V)$ and any proper normal subgroup of $C(V)$ is contained in $Z(V)$, i.e. $C(V)/ (C(V) \cap Z(V))$ is simple. The  proofs of  results announced by Dieudonn\'e are left to the reader.

We  prove the following generalization
\begin{thm}
The group $C(V)$ is simple. Moreover, if $D= \mathbb{F}_2$ is the two element field, then $Z(V)$ is trivial and $C(V)$ coincides with $F(V)$.
\end{thm}

For any infinite cardinal $\aleph_{\nu} \leq \aleph_{\delta}$,  by $F_{\nu}$ we denote the set of all linear transformations with ranges of dimension $<\aleph_{\nu}$. The known result says that all proper nontrivial  two-sided ideals of $End(V)$  have a form $F_{\nu}$ for some $0 \leq \nu  \leq \delta$ (see \cite{Jac}, Chapter IX). So, we have a chain of all two-sided ideals $\{0\} < F_0 <F_1 < \ldots < F_{\delta} < End(V)$ corresponding to cardinals $\aleph_0 < \aleph_1< \ldots < \aleph_{\delta}$. It is clear that this result can be generalized to the description of normal subgroups of 
$GL(V)$. It was done by A. Rosenberg in \cite{Ros58}. 

Let $Z^{*}$ be the multiplicative group of the center $Z$ of division ring $D$. 
 For any infinite cardinal $\aleph_{\nu} \leq \aleph_{\delta}$ by $G_{\nu(1)}$ we denote the  subgroup of $GL(V)$ consisting of all elements of the form $id_V + A$, where $ A \in F_{\nu}$. Clearly, $G_{\nu(1)}$ are normal subgroups of $GL(V)$. However, in contrast with the description of two-sided ideals in $End(V)$ the group of scalar operators also is a nontrivial  normal subgroup in $GL(V)$ (in case of $End(V)$ any two-sided ideal containing $id_V$ coincides with whole  $End(V)$). 
 Moreover, the group $G_{\nu} = Z^{*} \cdot id_V \times G_{\nu(1)}$ is normal in $GL(V)$.
 
 The main result (Theorem B) of \cite{Ros58} shows that: 
 \begin{quote} if $N$ is a normal subgroup of $GL(V)$, then either 
 
 (i) $N= H \cdot id_V \times G_{\nu(1)}$, where $0< \nu \leq \delta$ and $H \leq Z^{*}$, 
 
 or 
 
 (ii) $N \leq  Z^{*} \cdot id_V \times G_{0(1)} = G_0$.
 \end{quote}

This result describe completely the normal subgroups of $GL(V)$ which strictly contain $Z^{*} \cdot id_V \times G_{0(1)}$ (the case (i)). As noted Rosenberg, except of subgroups of $Z^{*} \cdot id_V \times G_{0(1)} = G_0$ the only normal subgroups that arise are due to the ideals of $End(V)$.  There is no doubt that the proof of case $(i)$ is correct. The case (ii) A. Rosenberg relates to Dieudonn\'e's paper,  he  wrote:
  \begin{quote}
  "Furthermore, in [5, p. 45] Dieudonn\'e studied the case (ii). He showed there, that either 
  $N \subseteq Z^{*} \cdot id_V$ or $N \supseteq
   (G_{0(1)} \cap C)= C_0$. Moreover, $C_0=G_{0(1)}'$, the group $C_0/ (C_0 \cap Z^{*} \cdot id_V)$ is simple,  and    $G_{0(1)}/ C_0$ is isomorphic to the multiplicative group of $D$ made abelian."
\end{quote} 
\noindent Here Rosenberg defines  $C_0=G_{0(1)} \cap C$, where $C$ is the group generated by elements of class two, i.e. elements  of the form $id_V + A$, where $A^2=0$ (see \cite{Baer}, page $207$),  and says (without proof) that $C_0$ is equal to commutator subgroup $G_{0(1)}'$ of $G_{0(1)}$. In Dieudonn\'e's notation  $G_{0(1)}=F(V)$ and $G_{0(1)}'= (F(V))'=C(V)$). 
 Theorem $A$ of \cite{Ros58} says that $GL(V)$ is generated by elements of class $2$. It means that $C = GL(V)$ and $C_0=G_{0(1)} \cap C=G_{0(1)}$. So, unfortunately the  statement  $C_0=G_{0(1)}'$  is  false and $C_0$ should be defined in other way.

The only new result of Rosenberg which extends,  in case $(ii)$, results of Dieudonn\'e is Lemma $3.7$ which says that every normal subgroup of
$G_{\nu} = Z^{*} \cdot id_V \times G_{\nu(1)}$ is normal in $GL(V)$. So it suffices to find only normal subgroups in  $Z^{*} \cdot id_V \times G_{0(1)}$ to complete description of normal subgroups of 
$GL(V)$.

If $V$ is countably dimensional over a field $K$ and we fix a basis $e_1, e_2, \ldots$ then $GL(V)$ corresponds to the group  $GL_{cf}(\mathbb{N}, K)$  of invertible column-finite infinite matrices over $K$ indexed by 
 $\mathbb{N}$.  
 Let $GL_{fr}(\mathbb{N}, K)$ denotes the subgroup of 
 $GL_{cf}(\mathbb{N}, K)$ which consits of all matrices which differ from the identity matrix only in finite numbers of rows. By $SL_{fr}(\mathbb{N}, K)$ we denote the subgroup fo $GL_{fr}(\mathbb{N}, K)$ consisting of all matrices for which determinant of the submatrix in left upper corner covering nonidentity rows is equal to $1$. 
By $D_{sc}(\mathbb{N}, K)$ we denote the subgroup of all scalar matrices in $GL_{cf}(\mathbb{N}, K)$.

It is clear that $GL_{fr}(\mathbb{N}, K)$, $SL_{fr}(\mathbb{N}, K)$ and $D_{sc}(\mathbb{N}, K)$   corresponds to $F(V) = G_{0(1)}$,  $C(V) = (G_{0(1)})'$ and  $Z(V)= Z^{*} \cdot id_V$ respectively (Rosenberg's and Dieudonn\'e's definitions are basis independent).

Rosenberg's  Theorem B says that all proper normal subgroups of $GL_{cf}(\mathbb{N}, K)$ are contained in $D_{sc}(\mathbb{N}, K) \times  GL_{fr}(\mathbb{N}, K)$. It means that 
 $D_{sc}(\mathbb{N}, K) \times  GL_{fr}(\mathbb{N}, K)$ is maximal normal subgroup of $GL_{cf}(\mathbb{N}, K)$ and the corresponding factorgroup is simple.  From Lemma 3.7 of \cite{Ros58} every normal subgroup of $D_{sc}(\mathbb{N}, K) \times  GL_{fr}(\mathbb{N}, K)$ is normal in $GL_{cf}(\mathbb{N}, K)$ (some kind of transitivity of normality).

Simplicity of $SL_{fr}(\mathbb{N}, K)$ was proved by Clowes and Hirsch in \cite{CloH} (not known to Rosenberg).  Moreover, this proof shows, as in finite dimensional case, that every normal subgroup of $GL_{cf}(\mathbb{N}, K)$ which is not contained in  $D_{sc}(\mathbb{N}, K)$ contains $SL_{fr}(\mathbb{N}, K)$. In the proof it was also shown that $SL_{fr}(\mathbb{N}, K)$ is generated by elementary transvections and its conjugates.

We use above results and results on generators of the group $GL_{cf}(\mathbb{N}, K)$ given in \cite{Hol-book}, \cite{Tho}, \cite{Ver} to prove the following:
\begin{thm}
The subgroups
\begin{quote}
 $D_{sc}(\mathbb{N}, K)$,\\ $SL_{fr}(\mathbb{N}, K)$,\\  $GL_{fr}(\mathbb{N}, K)$,\\  $D_{sc}(\mathbb{N}, K) \times   SL_{fr}(\mathbb{N}, K) $,\\  $D_{sc}(\mathbb{N}, K) \times  GL_{fr}(\mathbb{N}, K)$
\end{quote}
 are  normal subgroups of $GL_{cf}(\mathbb{N}, K)$.

The group  $SL_{fr}(\mathbb{N}, K)$ and the factor group $GL_{cf}(\mathbb{N}, K)/ (D_{sc}(\mathbb{N}, K) \times   GL_{fr}(\mathbb{N}, K)) $ are simple. The group $D_{sc}(\mathbb{N}, K)$ and the factor groups
\begin{quote} $ (D_{sc}(\mathbb{N}, K)  \times  SL_{fr}(\mathbb{N}, K)) / SL_{fr}(\mathbb{N}, K) $, \\
 $(D_{sc}(\mathbb{N}, K) \times   GL_{fr}(\mathbb{N}, K)) / GL_{fr}(\mathbb{N}, K) $,\\ $GL_{fr}(\mathbb{N}, K)/ SL_{fr}(\mathbb{N}, K)$, \\  $(D_{sc}(\mathbb{N}, K) \times   GL_{fr}(\mathbb{N}, K)) /  (D_{sc}(\mathbb{N}, K)  \times  SL_{fr}(\mathbb{N}, K))$
\end{quote}
  are isomorphic to $K^{*}$.
\end{thm}

 \begin{Cor}
$1)$ The group $GL_{cf}(\mathbb{N}, K)$ is perfect, i.e. it coincides with its commutator subgroup.

\noindent $2)$ $D_{sc}(\mathbb{N}, K)$ is the center of $GL_{cf}(\mathbb{N}, K)$.

\noindent $3)$ $SL_{fr}(\mathbb{N}, K)$ is the commutator subgroup of $GL_{fr}(\mathbb{N}, K)$.
\end{Cor}

In fact $GL_{cf}(\mathbb{N}, K)$ coincides with all terms of its lower central series.

Since $D_{sc}(\mathbb{N}, K)$ is the center of $GL_{cf}(\mathbb{N}, K)$ ( Lemma $2.4$), every subgroup of $D_{sc}(\mathbb{N}, K)$ is normal in $GL_{cf}(\mathbb{N}, K)$. Moreover, any subgroup $H$ such that $SL_{fr}(\mathbb{N}, K) \leq H \leq GL_{fr}(\mathbb{N}, K)$ is normal in $GL_{fr}(\mathbb{N}, K)$ (Lemma $2.5$).

The lattice of normal subgroups of $GL_{cf}(\mathbb{N}, K)$  ``modulo the center'' is shown in the  figure below (we abbreviate notation for convenience).
The thin line between subgroups $H_1$ and $H_2$ ($H_1 < H_2$) means that the factor group $H_2/ H_1$ is simple, the thick line means that the factor group $H_2/H_1$ is isomorphic to $ K^{*}$.    

\begin{center}
\begin{tikzpicture}
 
 \draw[ultra thick] (0,0)--(-2,1);
 \draw (-2,1)--(0,2.5);
 \draw[ultra thick] (0,2.5)--(2,4);
 \draw (2,4)--(2,5.5);
 \draw (0,0)--(2,1.5);
 \draw[ultra thick] (2,1.5)--(4,3);
 \draw[ultra thick] (4,3)--(2,4);
 \draw[ultra thick] (2,1.5)--(0,2.5);
 
 \fill (0,0) node [fill=white]{$\{ e \}$};
 \fill (-2,1) node [fill=white]{$D_{sc}$};
 \fill (2,1.5) node [fill=white]{$SL_{fr}$};
 \fill (4,3) node [fill=white]{$GL_{fr}$};
 \fill (0,2.5) node [fill=white]{$D_{sc} \times SL_{fr}$};
 \fill (2,4) node [fill=white]{$D_{sc} \times  GL_{fr}$};
 \fill (2,5.5) node [fill=white]{$GL_{cf}$};
 \fill (3,0.5) node [fill=white]{$\leftarrow$ {\rm Clowes, Hirsch}};
 \fill (4,4.7) node [fill=white]{$\leftarrow$ {\rm Rosenberg}};
 \fill (2,2.7) node [fill=white]{?};
\end{tikzpicture}
\end{center}
\noindent {\bf Remark 1.} If $K$ is a two element field $\mathbb{F}_2$, then $D_{sc}$ is trivial, $GL_{fr}=SL_{fr}$ and this lattice reduces  in obvious way.
$GL_{fr}(\mathbb{N}, \mathbb{F}_2)$ is the only proper nontrivial normal subgroup of $GL_{cf}(\mathbb{N}, \mathbb{F}_2)$. \qed

\smallskip

It is clear that Theorem $1.2$ reduces the problem of determining the normal subgroups of  $GL_{cf}(\mathbb{N}, K)$ to the following problem

\begin{quote}
{\bf Problem} Describe all normal subgroups of the group $D_{sc}(\mathbb{N}, K) \times GL_{fr}(\mathbb{N}, K)$ which contain $SL_{fr}(\mathbb{N}, K)$.
\end{quote}

The partial answer is given by
\begin{prop}
If $H_1 \leq D_{sc}(\mathbb{N}, K)$ and $H_2$ is such that 
$SL_{fr}(\mathbb{N}, K) \leq H_2 \leq GL_{fr}(\mathbb{N}, K)$, then 
$H_1 \times H_2$ is a normal subgroup of $D_{sc}(\mathbb{N}, K) \times GL_{fr}(\mathbb{N}, K)$.

\end{prop}

\noindent {\bf Remark 2.} Rosenberg omitted many normal subgroups described in Proposition $1.4$. He only noted at the beginning of paragraph $3$ of \cite{Ros58} that $H_1 \times GL_{fr}(\mathbb{N}, K)$ are normal in $GL_{cf}(\mathbb{N}, K)$ if $H_1 \leq D_{sc}(\mathbb{N}, K)$) but he omitted all subgroups of the form $H_1 \times H$, where $SL_{fr}(\mathbb{N}, K) \leq H  <  GL_{fr}(\mathbb{N}, K)$ (with exception of $SL_{fr}(\mathbb{N}, K)$ of course). \qed

\smallskip

Moreover, the following example shows that there exist solutions of our problem which are not a direct product.

\smallskip

\noindent {\bf Example} From Theorem $1.2$ it follows  $D_{sc}(\mathbb{N}, K) \cong K^{*}$. 
Let $f$ be a natural homomorphism $f: GL_{fr}(\mathbb{N}, K) \to 
GL_{fr}(\mathbb{N}, K)/SL_{fr}(\mathbb{N}, K) \cong K^{*}$. 
Let 
$$H = \{ (a,b) \in  D_{sc}(\mathbb{N}, K) \times GL_{fr}(\mathbb{N}, K): f(b)=a \}.$$
One can easily check that $H$ is normal subgroup in $D_{sc}(\mathbb{N}, K) \times GL_{fr}(\mathbb{N}, K)$. Since $H \cap D_{sc}(\mathbb{N}, K)$ is trivial and $H \cap SL_{fr}(\mathbb{N}, K) =SL_{fr}(\mathbb{N}, K)$ we deduce that $H$ is not a direct product.
In fact $H$ is a subdirect  product. \qed

\smallskip

The main idea in solving the problem uses the fact that if $H$ is subgroup of abelian group $G$ and $\varphi: H \to \varphi (H) \subset G$ is an isomomorphism, then 
$\Delta_{\varphi}(H) = \{ (g, \varphi(g)): g \in H\}$ is a subgroup of $G \times G$. 

By $\pi_i$ we denote the natural projection onto $i-th$ component of direct product. We give now a complete answer to our problem.

\smallskip

\begin{thm}
If $H$ is a normal subgroup of $D_{sc}(\mathbb{N}, K) \times  GL_{fr}(\mathbb{N}, K)$ containing  $SL_{fr}(\mathbb{N}, K)$ then $H$ is determined in a unique way by a quintuple $(H_1, F_1, H_2, F_2, \varphi)$ where $H_1 = H \cap D_{sc}(\mathbb{N}, K)$, $H_2 = H \cap  GL_{fr}(\mathbb{N}, K)$, $F_i = \pi_i (H)$, $\varphi: F_1/ H_1 \to F_2/H_2$ is an isomorphism 
and $(x,y) \in H$ if and only if  $\varphi(xH_1)=yH_2$.
\end{thm}

It is clear that the quintuple $(H_1, H_1, H_2, H_2, \varphi)$ corresponds to $H= H_1 \times H_2$. In general $H$ is only a subdirect product of $\pi_1(H) \times \pi_2(H)$ (which is a fibre product by Goursat's Lemma \cite{G}).


If $V$ is countably dimensional over a field $K$ it is clear that Rosenberg omitted many normal subgroups described in Proposition $1.4$ which are direct product and all normal subgroups which are not direct product (see Example and Theorem $1.5$). 
Moreover, Rosenberg and Dieudonn\'e did not remarked the  triviality of normal subgroup $D_{sc}(\mathbb{N}, K) \cap  SL_{fr}(\mathbb{N}, K)$ and simplification of description in case $K=\mathbb{F}_2$ (Remark $1$, Theorem $1.1$).   
So,   Rosenberg's results are  incomplete. 
 
Our Theorems $1.2$ and $1.5$ fill these gaps  giving the full  description of the lattice of normal subgroups of the group of infinite column-finite matrices indexed by positive integers over any field.

 We believe that for $V$ which is not countably dimensional the description of normal subgroups of $GL(V)$ contained in $Z^{*} \cdot id_V \times G_{0(1)}$ is similar but it need more technical work (there are no references like \cite{BH}, \cite{CloH}, \cite{Hol-book}, \cite{Tho}, \cite{Ver} for countable case over a field). Since every finite division ring is a field by Wedderburn's Theorem, our resuts hold also for finite division rings.

We note that similar description of ideals of the Lie algebra of infinite column-finite matrices over any field was obtained in \cite{Hol-inf-Lie}.

\section{Proofs of main results}

\noindent {\it Proof of Theorem $1.1$.}

Assume that $C(V)/ (C(V) \cap Z(V))$ is simple. It is clear that the center  $Z(V)=Z^{*} \cdot id_V$ of $GL(V)$ contains only homotheties and for every nonidentity  homothety $g$ the element $g-id_V$  has infinite dimensional range if 
 $V$ is infinite dimensional  and $D \not= \mathbb{F}_2$.  If $D = \mathbb{F}_2$ the center is trivial. So $C(V) \cap Z(V)$ is  always trivial  and $C(V)$ is simple for all division rings. \qed

Let $K$ be any  field. We start with the following Lemmas

\begin{lem}
If the subgroup $A_1$ is normal in the group $A$ and the subgroup $B_1$ is normal in the group $B$, then $A_1 \times B_1$ is normal in $A \times B$  and $$(A \times B)/ (A_1 \times B_1) \cong (A/A_1) \times (B/B_1).$$
\end{lem}

The above Lemma  is a Theorem $2.30$ from \cite{Rotman}.

\begin{lem}
The group  $GL_{cf}(\mathbb{N}, K)$ is generated by the row and column-finite matrices and  upper triangular matrices.
\end{lem}

This lemma follows easily from considerations in chapter $2$ of \cite{Ros58}. Another proof one can find in \cite{Ver}.

An infinite block diagonal matrix with finite blocks of sizes $n_1 \times n_1$,  $n_2 \times n_2$, $\ldots$ is called a {\it string} with the shape $(n_1, n_2, \ldots)$. Of course, the string with shape $(1,1, \ldots)$ is a diagonal matrix. 
The following Lemma is Theorem $3.3$ of \cite{Tho}. Its proof is contained in chapter $5$ of \cite{Tho}. P. Vermes in \cite{Ver} gave another proof for $K=\mathbb{C}$, however it can be easily adopted to arbitrary field $K$ (see also \cite{Hol-book}). By $GL_{rcf}(\mathbb{N}, K)$ we denote the subgroup of $GL_{cf}(\mathbb{N}, K)$ consisting of the row and column-finite matrices.
\begin{lem}
The group  $GL_{rcf}(\mathbb{N}, K)$ is generated by strings.
\end{lem}

\smallskip

\noindent {\it Proof of Theorem $1.2$.}

Now we prove normality of groups listed in Theorem $1$. We note that by Lemma $3.7$ of \cite{Ros58} it suffices to prove normality in $D_{sc}(\mathbb{N}, K) \times  GL_{fr}(\mathbb{N}, K)$.
The normality of $D_{sc}(\mathbb{N}, K)$ follows from Lemma $2.1$. Let $ g \in SL_{fr}(\mathbb{N}, K)$ and $s$ be a string with shape $(n_1, n_2, \ldots)$. Then $g$ has a form 
$$\left(\begin{array}{c|c} \hat{g}& \star \cr \cline{1-2}
0&e\cr
\end{array}\right)$$
where $\hat{g} \in SL(n, K)$ for some $n \in \mathbb{N}$ and $e$ is the infinite identity matrix. We choose minimal $t$ such that $m = \sum_{i=1}^{t} n_i \geq n$.
We can extend block decomposition of $g$ to 

$$\left(\begin{array}{c|c} \tilde{g}& \star \cr \cline{1-2}
0&e\cr
\end{array}\right)$$
where $\tilde{g} \in SL(m, K)$ .
Then $s^{-1} gs$ has a form
$$\left(\begin{array}{c|c} \hat{g} & \star \cr \cline{1-2}
0&e\cr
\end{array}\right)$$
where $\hat{g} \in SL(m,K)$.
If $u$ is any upper triangular matrix, then we can use on $u$ the same block structure as on $g$. 
Simple calculations show that  $u^{-1} gu \in SL_{fr}(\mathbb{N}, K)$.
Now, if $h=h_1 \cdot h_2 \cdot  \ldots \cdot h_n$ is the product of upper triangular and block-diagonal matrices, then we deduce 
$h^{-1} gh \in SL_{fr}(\mathbb{N}, K)$ by induction on $n$.
So $SL_{fr}(\mathbb{N}, K)$ is normal in $GL_{cf}(\mathbb{N}, K)$ and of course in $GL_{fr}(\mathbb{N}, K)$.

 The normality of other subgroups from Theorem $1.2$ follows from Lemma $2.1$. 

Simplicity of $SL_{fr}(\mathbb{N}, K)$ was proved by Clowes and Hirsch in \cite{CloH}.  Moreover, this proof shows, as in finite dimensional case, that every normal subgroup of $GL_{cf}(\mathbb{N}, K)$ which is not contained in  $D_{sc}(\mathbb{N}, K)$ contains $SL_{fr}(\mathbb{N}, K)$.
The fact that the factor group 
$GL_{cf}(\mathbb{N}, K)/ (D_{sc}(\mathbb{N}, K) \times   GL_{fr}(\mathbb{N}, K)) $ 
is simple follows from Theorem B of \cite{Ros58}. 

Let $d(\alpha)= (\alpha -1) e_{11} + e$, where $e$ has $1$ in $(1,1)$ entry and zero elsewhere. 
For every matrix $g \in GL_{fr}(\mathbb{N}, K)$ we have a unique decomposition
$$
g=d(\alpha) \cdot (d(\alpha^{-1}) \cdot g),
$$
where $\alpha = det(\hat{g})$ and $d(\alpha^{-1}) \cdot g \in SL_{fr}(\mathbb{N}, K)$. This shows that
 $GL_{fr}(\mathbb{N}, K)/ SL_{fr}(\mathbb{N}, K)$ is isomorphic to $K^{*}$. 
Now other statements of Theorem $1.2$ are obvious. \qed

\begin{lem}\label{center}
$D_{sc}(\mathbb{N}, K)$ is the center of $GL_{cf}(\mathbb{N}, K)$.
\end{lem}
As in a case of finite dimensional matrices we  first show that nondiagonal entries of the matrix from the center are equal to $0$, next we show that diagonal entries are equal. We  will omit the details.

\begin{lem}\label{norsl}
A nontrivial subgroup $H$ of $GL_{fr}(\mathbb{N}, K)$ is normal if and only if it  contains $SL_{fr}(\mathbb{N}, K)$.
\end{lem}

\noindent {\it Proof of Lemma $2.5.$}   

If $H$ is nontrivial normal subgroup of $GL_{fr}(\mathbb{N}, K)$, then from \cite{CloH} it follows that $H \geq SL_{fr}(\mathbb{N}, K)$. Since $GL_{fr}(\mathbb{N}, K) / SL_{fr}(\mathbb{N}, K) \cong K^{*}$ is abelian, $(GL_{fr}(\mathbb{N}, K))' = [GL_{fr}(\mathbb{N}, K), GL_{fr}(\mathbb{N}, K)] \leq SL_{fr}(\mathbb{N}, K)$ holds. If $SL_{fr}(\mathbb{N}, K) \leq H \leq GL_{fr}(\mathbb{N}, K)$, then 
$[GL_{fr}(\mathbb{N}, K), H] \leq [ GL_{fr}(\mathbb{N}, K), GL_{fr}(\mathbb{N}, K)] \leq  SL_{fr}(\mathbb{N}, K) \leq H$. So, $H$ is normal in $ GL_{fr}(\mathbb{N}, K)$. \qed

\bigskip

\noindent {\it Proof of Corollary $1.3.$}

1) In \cite{BH} it was shown that the uppertriangular matrix which has $1-$s on the main diagonal and on the second diagonal above the main diagonal and $0-$s otherwise belongs to the commutator subgroup of the group of upper unitriangular matrices  
$UT( \mathbb{N}, K)$. It does not belongs to $GL_{fr}(\mathbb{N}, K)$.
Since $(UT( \mathbb{N}, K))' \leq (GL_{fr}(\mathbb{N}, K))'$ the group 
 $GL_{cf}(\mathbb{N}, K)$ is perfect.

 2) follows from Lemma $2.4.$
 
3)  Simple calculations shows that  $(GL_{fr}(\mathbb{N}, K))' \leq SL_{fr}(\mathbb{N}, K)$. The reverse inclusion follows from Lemma $2.5.$
\qed

\bigskip

\noindent {\it Proof of Theorem $1.5$.}

Let $H$ be a normal subgroup of $D_{sc}(\mathbb{N}, K) \times  GL_{fr}(\mathbb{N}, K)$ containing  $SL_{fr}(\mathbb{N}, K)$. So $\pi_1(H)$ is normal in $D_{sc}(\mathbb{N}, K)$, $\pi_2(H)$ is normal in $GL_{fr}(\mathbb{N}, K)$ and $H$ is a subdirect product of $\pi_1(H) \times \pi_2(H)$.
Subgroups $\pi_2(H)$ are in $1-1$ correspondence with subgroups of $K^{*}$ by canonical homomorphism $GL_{fr}(\mathbb{N}, K) \to GL_{fr}(\mathbb{N}, K)/SL_{fr}(\mathbb{N}, K)$.

So, our problem is equivalent to the problem of description of all subgroups $H$ of abelian group $A_1 \times A_2$, where $A_1 \cong A_2 \cong K^{*}$. Let $H_i = H \cap A_i$. Since $H_1 \times H_2 \leq H$ we can factorize $(A_1 \times A_2)/ (H_1 \times H_2)$ and reduce by Lemma $2.1$ the problem to finding subgroups $F \leq (A_1/H_1) \times (A_2/H_2)$ such that
$F_i = F \cap (A_i/H_i) \ \ \ {\rm is \ \  trivial}.$
Let $K_i= \pi_i(F)$. 

 From definitions of $F_i$ it follows that for every $x \in K_1$ there exists a unique $y \in K_2$ such that $(x,y) \in F$. Moreover, for every $y \in K_2$ there exists a unique $x \in K_1$ such that $(x,y) \in F$. So, the mapping $x \mapsto y$ is an isomorphism of $K_1$ with $K_2$.

Now assume that we have $H_i \leq F_i \leq A_i$ and isomorphism $\varphi: F_1/H_1 \to F_2/H_2$. Then there exists a unique subgroup $H$ such that $H_i = H \cap A_i$, $F_i = \pi(H)$ and $(x,y) \in H$ if and only if $\varphi(xH_1)=yH_2$. 

Applying the above considerations in our case we finish the proof. \qed

{\bf Acknowledgements}

We are very gratful to A. V. Stepanov for a fruitful discussion and to anonymous referee for valuable remarks.


\end{document}